\documentclass[12pt]{article}
\usepackage{amsfonts, amsmath, latexsym, vmargin, epic, eepic}
\usepackage{epsfig}
\title{{\bf On simplicial and cubical complexes with short links}}
\setpapersize{custom}{21cm}{29.7cm}
\setmarginsrb{2cm}{1.5cm}{2cm}{3.5cm}{0pt}{0pt}{0pt}{0pt}

\author{Michel DEZA\\
CNRS/ENS, Paris and Institute of Statistical Mathematics, Tokyo,\\
\ Mathieu DUTOUR \\
ENS, Paris and Hebrew University, Jerusalem,
\footnote{Research of the second author was financed by EC's IHRP Programme, within the Research Training Network ``Algebraic Combinatorics in Europe,'' grant HPRN-CT-2001-00272.}\\
\ Mikhail SHTOGRIN 
\thanks{Third author acknowledges financial support of the Russian Foundation of Fundamental Research (grant 02-01-00803) and the Russian Foundation for Scientific Schools (grant NSh.2185.2003.1)}\\
Steklov Mathematical Institute, Moscow, Russia.}

\begin{document}
\newtheorem{theorem}{Theorem}[section]
\newtheorem{definition}[theorem]{Definition}
\newtheorem{conjecture}[theorem]{Conjecture}
\newtheorem{corollary}[theorem]{Corollary}
\newtheorem{problem}[theorem]{Problem}
\newtheorem{lemma}[theorem]{Lemma}
\newtheorem{claim}[theorem]{Claim}
\newtheorem{remark}[theorem]{Remark}
\newtheorem{proposition}[theorem]{Proposition}
\newcommand{\R}{\ensuremath{\mathbb{R}}}
\newcommand{\N}{\ensuremath{\mathbb{N}}}
\newcommand{\Q}{\ensuremath{\mathbb{Q}}}
\newcommand{\C}{\ensuremath{\mathbb{C}}}
\newcommand{\Z}{\ensuremath{\mathbb{Z}}}
\newcommand{\T}{\ensuremath{\mathbb{T}}}
\newcommand{\K}{\ensuremath{\mathbb{K}}}
\newcommand{\qed}{\hfill $\Box$ }
\newcommand{\proof}{\noindent{\bf Proof.}\ \ }

\maketitle

\begin{abstract}
We consider closed simplicial and cubical $n$-complexes in terms of link of their $(n-2)$-faces. Especially, we consider the case, when this link has size $3$ or $4$, i.e., every $(n-2)$-face is contained in $3$ or $4$ $n$-faces. 
Such simplicial complexes with {\em short} (i.e. of length $3$ or $4$) links are completely classified by their {\em characteristic partition}.
We consider also embedding into hypercubes of the skeletons of simplicial and cubical complexes.
\end{abstract}

\section{Introduction}
A {\em $n$-dimensional simplicial complex} (or {\em simplicial $n$-complex}) is a collection ${\cal S}$ of finite nonempty sets, such that:

(i) if $S$ is an element of ${\cal S}$, then so is every nonempty subset of $S$;

(ii) If $S,S'\in {\cal S}$, then $S\cap S'\in {\cal S}$;

(iii) all maximal (for inclusion) elements of ${\cal S}$ have cardinality $n+1$.

Given a simplicial complex ${\cal K}$ of dimension $n$, every $(k+1)$-subset of it defines a {\em face of dimension $k$} of the simplicial complex. In the sequel we identify faces with set of its vertices.
A simplicial complex ${\cal K}$ is called a {\em pseudomanifold} if every $(n-1)$-face belongs to one or two $n$-faces. The {\em boundary} is the set of $(n-1)$-faces contained in exactly one $n$-face. If the boundary of ${\cal K}$ is empty (i.e., every $(n-1)$-face is the intersection of exactly two $n$-faces), then ${\cal K}$ is called {\em closed} pseudomanifold.
%A simplicial complex is {\em closed} if every $(n-1)$-face is contained in exactly two $n$-faces.

For every $n$-face $F'$, containing $(n-2)$-face $F$, there exists unique edge $e$, such that $F'=F\cup e$. The {\em link of an $(n-2)$-face} $F$ is the family of cycles formed by above edges $e$, where $F'$ run through all $n$-simplexes, containing $F$.
If the number of such faces $F'$ is at most $5$ or if $K$ is a manifold, then link consists of unique cycle.
The length of this cycle will be denoted by $l(F)$.

Every compact manifold $M$ can be represented as a closed simplicial complex, if one choose a triangulation of $M$.

%NEED TO CONSIDER GENERAL REFERENCE FOR SIMPLICIAL COMPLEX

\begin{definition}
A closed simplicial complex ${\cal K}$ is called {\em of type $L$}, if $L$ is the set of all values $l(F)$, where $F$ is any $(n-2)$-face of ${\cal K}$.
\end{definition}
We will be concerned below, especially, with the case $L=\{3,4\}$. 

Some examples:

(i) the $n$-simplex (respectively, the $n$-hyperoctahedron) are examples of simplicial complexes of type $\{3,4\}$, where, moreover, $L=\{3\}$ (respectively, $L=\{4\}$);

(ii) the only $2$-dimensional $\{3,4\}$-simplicial complexes are: dual Triangular Prism, Tetrahedron and Octahedron.

Take the hyperoctahedron (also called {\em cross-polytope}) of dimension $n$ and write its set of vertices as $\{1,2,\dots, n+1, 1', 2', \dots, (n+1)'\}$. This $n$-hyperoctahedron has the following $2^{n+1}$ facets:
\begin{equation*}
\{x_1, \dots, x_{n+1}\}\mbox{~~with~~}x_i=i\mbox{~or~}i'\mbox{~for~}1\leq i\leq n+1\;.
\end{equation*}

\begin{definition}
Let $P=(P_1, \dots, P_t)$ be a partition of $V_{n+1}=\{1,2,\dots, n+1\}$.
Define the simplicial complex ${\cal K}(P)$ as follows:

(i) It has $n+1+t$ vertices $(1, 2,\dots, n+1, P_1,\dots,P_t)$.

(i) Every face $F_x=\{x_1,\dots x_{n+1}\}$ is mapped onto $F_y=\{y_1, \dots, y_{n+1}\}$, where $y_i=i$ if $x_i=i$, and $y_i=P_j$ if $x_i=i'$, $i\in P_j$. If two elements of $F_y$ are identical, then we reject the set $F_y$; otherwise, we add it to the list of facets of ${\cal K}(P)$.
\end{definition}

Below $K_m$ denotes the complete graph on $m$ vertices, $C_m$ denotes the cycle on $m$ vertices.
Denote by $K_m-C_h$ the complement in $K_m$ of the cycle $C_h$;
denote by $K_m-hK_2$ the complete graph on $m$ vertices with $h$ disjoint edges deleted.

Any simplicial or cubical complex of type $\{3,4\}$ is realizable as a manifold, since the neighborhood of every point is homeomorphic to the sphere.

Given a complex ${\cal K}$, the skeleton $G({\cal K})$ is the graph of vertices of ${\cal K}$ with two vertices adjacent if they form a $1$-face of ${\cal K}$.

Given a graph $G$, the {\em path-metric} (denoted by $d_G(i,j)$) between two vertices $i$, $j$ is the length of a shortest path between them.
The graph $G$ is said to be {\em embeddable up to scale $\lambda$ into hypercube} if there exist a mapping $\phi$ of $G$ into $\{0,1\}^N$ with $|\phi(i)-\phi(j)|=\lambda d_G(i,j)$.
For the details on such embeddability, see the book \cite{DL}.

For example, Proposition 7.4.3 of \cite{DL} give that 
$K_{m+1}-K_2$ and $K_{2m}-mK_2$ embed in the $2a_m$-hypercube with a scale $\lambda=a_m$, where $a_m={m-2\choose \frac{m}{2}-1}$ for $m$ even and $a_m=2{m-2\choose \frac{m-3}{2}}$ for $m$ odd.
Clearly, any subgraph $G$ of $K_{2n}-nK_2$, containing $K_{n+1}-K_2$, also admits above embedding, since any subgraph of diameter two graph is an {\em isometric} subgraph. 
In general, if $G$ is isometric subgraph of an hypercube, then it is an induced subgraph but this implication is strict.

A graph is said to be {\em hypermetric} if its path-metric satisfies the inequality
\begin{equation*}
\sum_{1\leq i<j\leq n} b_ib_j d_{G}(i,j)\leq 0
\end{equation*}
for any vector $b\in \Z^n$ with $\sum_{i}b_i=1$. In the special case, when $b$ is a permutation of $(1,1,1,-1,-1,0,\dots,0)$, above inequality is called {\em $5$-gonal}. The validity of hypermetric inequalities is necessary for embeddability but not sufficient: an example of hypermetric, but not embeddable graph (amongst those, given in Chapter 17 of \cite{DL}) is $K_7-C_5$.

\section{Simplicial complexes of type $\{3,4\}$}

In this section, we classify the simplicial complexes
of type $\{3,4\}$ in terms of partitions. Let $K$ be a simplicial complex
of type $\{3,4\}$ and let $\Delta=\{1,\dots,n+1\}$ be a $n$-face
of this complex. Denote by $F_{i}=\{1,\dots, i-1, i+1,\dots,n+1\}=V_{n+1}-\{i\}$
a facet of $\Delta$. This facet is contained in another simplex,
which we write as $\Delta_i=\{1,\dots,i-1,i',i+1, \dots, n+1\}$.
Denote by $F_{i,j}=V_{n+1}-\{i,j\}$ the $(n-2)$-faces of $K$.
One has $l(F_{i,j})=3$ if and only if $i'=j'$. Now, $l(F_{i,j})=4$
if and only if $i'\not= j'$ and $(i',j')$ is an edge.

%The {\em characteristic} of the simplex $\Delta$ is the expression
%$\phi(\Delta)=(l(F_{i,j}))_{1\leq i<j\leq n+1}$.
Define a graph on the set $V_{n+1}$ by making $i$ and $j$
adjacent if $l(F_{i,j})=3$. By what we already know
about $i'$ and $j'$, one obtains that this graph is of the form
$K_{P_1}+\dots+K_{P_t}$ (where $K_{A}$ denotes the complete graph on the
vertex-set $A$) and so, one gets a {\em characteristic partition}
of $V_{n+1}$, which we write as $P=\{P_1, \dots, P_t\}$.

%NEED TO DO A RETHINKING OF THE PROOF

Now we want to prove the equality between $K$ and the polytope ${\cal K}(P)$. In fact, it follows from the following lemma.

\begin{lemma}\label{FundamentalLemma}
Given a simplex $\Delta=\{1,2,\dots, n+1\}$ and a characteristic partition
of $\Delta$, there exist an unique complex $K$ of type $\{3,4\}$ having
$\Delta$ as a simplex.

Moreover, all simplexes of $K$ have the same characteristic partition.
\end{lemma}
\proof According to above notations, we define the vertices $i'$ and simplexes $\Delta_{i}$, such that $\Delta\cap \Delta_{i}=F_i$. The vertex-set of the complex $K$ contains vertices $\{1, \dots, n+1, 1', \dots, (n+1)'\}$; we will show below that it contains no others.

Let us find the values of the numbers $l(F)$ for the simplex, say, $\Delta_1$.

Take a $(n-2)$-face $F$ in $\{1',2,\dots, n+1\}$. If $1'\notin F$, then one has a $(n-2)$ face of $\Delta$ and so, we already know $l(F)$.

Let us write $F$ as $F'_{i,j}=\{1',2,\dots, n+1\}-\{i,j\}$. 
The face $F_{i,j}$ is contained in the simplexes $\Delta$,
$\Delta_i$ and $\Delta_{j}$.
If $l(F_{i,j})=3$, then $i'=j'$.
If $l(F_{i,j})=4$, then $F_{i,j}$ is also contained in $F_{i,j}\cup \{i', j'\}$, which is a simplex of $K$.

The face $F'_{i,j}$ is contained in the $(n-1)$-faces
$\{1',2,\dots, n+1\}-\{i\}$ and $\{1',2,\dots, n+1\}-\{j\}$.
According to $l(F_{1,i})=3$ or $4$, the face $F'_{i,j}$
is contained in either $\Delta_i$ (and $i'=1'$), or in
$\{1',2,\dots, i-1, i', i+1, \dots, n+1\}$. The same holds for $F_{1,j}$.

We deal with all cases:
\begin{itemize}
\item If $l(F_{i,j})=3$, one has $i'=j'$.
\begin{itemize}
\item If $l(F_{1,i})=4$ and $l(F_{1,j})=4$, then $F'_{i,j}$ is contained in $\{1',2,\dots, n+1\}$, $\{1',2,\dots, i-1, i', i+1, \dots, n+1\}$ and $\{1',2,\dots, j-1, j', j+1, \dots, n+1\}$. 
By equality $i'=j'$, one has $l(F'_{i,j})=3$.
\item If $l(F_{1,i})=3$, then $1'=i'$; so, $1'=j'$ and $l(F_{1,j})=3$. The face $F'_{i,j}$ is contained in $\{1',2,\dots, n+1\}=F'_{i,j}\cup \{i,j\}$, 
$\{1',2,\dots, i-1, 1, i+1, \dots, n+1\}=F'_{i,j}\cup \{1,j\}$
and $\{1',2,\dots, j-1, 1, j+1, \dots, n+1\}=F'_{i,j}\cup \{i,1\}$.
By equality $i'=j'$, one has $l(F'_{i,j})=3$.
\end{itemize}
\item If $l(F_{i,j})=4$, one has $i'\not=j'$.
\begin{itemize}
\item If $l(F_{1, i})=4$ and $l(F_{1,j})=4$, then $F'_{i,j}$ is contained in $\{1',2,\dots, n+1\}=F'_{i,j}\cup \{i,j\}$, $\{1',2,\dots, i-1, i', i+1,\dots, n+1\}=F'_{i,j}\cup \{i',j\}$, $\{1', 2, \dots, j-1, j',j+1, \dots, n+1\}=F'_{i,j}\cup \{i,j'\}$. Since the length of link should be $3$ or $4$ and we have already $4$ vertices, one gets that $F'_{i,j}$ is contained in $F'_{i,j}\cup \{i',j'\}$ and $l(F'_{i,j})=4$.
\item If $l(F_{1,i})=3$, then $1'=i'$ and so, $1'\not= j'$, which implies $l(F_{1,j})=4$. 
The face $F'_{i,j}$ is contained in $\{1',2,\dots, n+1\}=F'_{i,j}\cup\{i,j\}$, 
$\{1',\dots, i-1, 1,i+1,\dots, n+1\}=F'_{i,j}\cup \{1,j\}$ and $\{1',2,\dots, j-1, j',j+1, \dots, n+1\}=F'_{i,j}\cup \{i,j'\}$. So, by the same argument, one gets that $F'_{i,j}$ contained in $F'_{i,j}\cup \{1,j'\}$ and $l(F'_{i,j})=4$.
\end{itemize}
\end{itemize}
One obtains $l(F'_{i,j})=l(F_{i,j})$. Therefore, $\Delta$ and $\Delta_1$ have the same characteristic partition. Moreover, one see that the adjacent simplexes to $\Delta_1$ are contained in the vertex-set ${\cal V}=\{1, \dots, n+1, 1',\dots, (n+1)'\}$. This implies that the vertex-set of the defined complex is exactly ${\cal V}$. So, the complex of type $\{3,4\}$ $K$ is uniquely defined by the characteristic partition of any one of its simplexes. \qed

%types $\{3,4\}$ of the facet $\Delta_1$ are the same as the type $\{3,4\}$ of the facet $\Delta$. In particular, one can compute those types and so, the complex $K$, from the types of the simplex $\Delta$. \qed

Given a complex ${\cal K}$, its automorphism group $Aut({\cal K})$ is defined as the group of permutations of its vertices, preserving the set of faces.

Call a complex {\em isohedral} if $Aut({\cal K})$ is transitive on its facets.

\begin{corollary}\label{ConsequenceFundamentalLemma}
(i) All simplicial complexes of type $\{3,4\}$ are of the form ${\cal K}(P)$;

(ii) given two partitions $P$ and $P'$ of $V_{n+1}$, one has ${\cal K}(P)$ isomorphic to ${\cal K}(P')$ if and only if $P'$ is obtained from $P$ by a permutation of $V_{n+1}$;

(iii) every simplicial complex of type $\{3,4\}$ is isohedral.
\end{corollary}
\proof (i) Take a simplicial complex $K$ of type $\{3,4\}$ and a simplex $\Delta=\{1, 2\dots, n+1\}$ in it. This simplex has a characteristic partition $P$ and so, coincide with ${\cal K}(P)$.

(ii) By Lemma \ref{FundamentalLemma}, all simplexes of a simplicial complex of type $\{3,4\}$ have the same characteristic partition. So, if two complexes of type $\{3,4\}$ are isomorphic, their corresponding partitions are isomorphic too. 
On the other hand, two isomorphic partitions define the same simplicial complex of type $\{3,4\}$.

(iii) By (i), one can assume that $K$ is of the form ${\cal K}(P)$.
Take another facet $\Delta'=\{v_1,\dots, v_{n+1}\}$ of ${\cal K}(P)$, its partition type is the same as of $\{1,2,\dots, n+1\}$. So, one can construct a mapping $\varphi$ from $\{1, \dots, n+1\}$ to $\{v_1, \dots, v_{n+1}\}$ preserving partitions and, by extension being an automorphism of the complex. \qed

One can check, that if ${\cal K}$ is a simplicial complex of type $\{3,4\}$, such that $l(F)=3$ for each $(n-2)$-faces $l(F)$, contained in a fixed $(n-1)$-face, then ${\cal K}$ is the simplex.

\begin{theorem}
Every simplicial complex of type $\{3,4\}$ is spherical.
\end{theorem}
\proof Take the trivial partition $\{1\},\dots, \{n+1\}$ of $V_{n+1}$, the corresponding complex is hyperoctahedron, which is spherical, of course. Take now a complex ${\cal K}(P)$ of type $\{3,4\}$. By merging vertices $i'$ and $j'$, belonging to the same part, we preserve the sphericity. Furthermore, while doing this operation, we do not obtain pair of different faces having the same set of vertices. So, obtained simplicial complexes are necessarily spherical. \qed

\begin{proposition}
The skeleton of the simplicial $n$-complex ${\cal K}(P)$ of type $\{3,4\}$ is $K_{n+1+t}-hK_2$, where $t$ is the number of sets of the partition and $h$ is the number of singletons in the partition.
\end{proposition}
\proof In the hyperoctahedron, a point $i$ is not adjacent only to the point $i'$. If $i'$ belongs to a partition of size bigger than $1$, then an edge appears; otherwise, there is no such edge. \qed

In particular, $K_{n+1}$ and $K_{2n}-nK_2$ are the skeletons of, respectively, $n$-simplex and $n$-hyperoctahedron.

\begin{theorem}\label{ProductDecomposition}
Every simplicial polytope of the form ${\cal K}(P)$ with $P=(P_1, \dots, P_t)$ is isomorphic to the dual of the product complex $\Delta_1\times\dots\times\Delta_t$ with $\Delta_i$ being the simplex of dimension $|P_i|$.
\end{theorem}
\proof Let us denote $K=\Delta_1\times\dots\times\Delta_t$; by Corollary \ref{ConsequenceFundamentalLemma}, it suffices to prove that the dual complex $K^*$ is of type $\{3,4\}$ and that its characteristic partition is $P$.

We are reasoning in dual terms; we fix a vertex $1_i\in \Delta_i$ and so, a vertex $v=(1_1, \dots, 1_t)\in K$. Any two adjacent vertices in $K$ differ by exactly one coordinate. Take $F$ a $2$-dimensional face of $K$, which contains $v$; the vertex $v$ is adjacent to two vertices $v_1$ and $v_2$ contained in $F$. Denote by $x_i$ the coordinates of $v_i$, which differ from $v$. If $x_1=x_2$, then $v_1$ is adjacent to $v_2$ and so, $F$ has three vertices. If $x_1\not= x_2$, then $F$ contains the vertex $v'=(1_1,\dots, (v_1)_{x_1},\dots, (v_2)_{x_2}, \dots, 1_t)$, where $(v_i)_{x_i}$ denotes the $x_i$-th coordinate of $v_i$. So, $F$ contains the four vertices $v$, $v_2$, $v'$, $v_1$, which form a square. Therefore, $K^*$ is of type $\{3,4\}$.

On the other hand, above computation proves that $K^*$ is isomorphic to ${\cal K}(P)$. \qed

\begin{proposition}
Every simplicial complex ${\cal K}$ of type $\{3,4\}$ admits a polytopal realization, such that its group of isometries coincides with its group $Aut({\cal K})$ of combinatorial isometries.
\end{proposition}
\proof By theorem \ref{ConsequenceFundamentalLemma}, one can assume that $K={\cal K}(P)$. This result follows immediately from the product decomposition given in Proposition \ref{ProductDecomposition}. \qed

Note that, in general, this group of isometries is a subgroup of $Aut({\cal K})$. Above proposition is an analog of Mani's theorem (\cite{M}) for $3$-connected planar graphs.

\begin{proposition}
Take a partition $P=(P_1,\dots, P_t)$ of $V_{n+1}$; the order of $Aut({\cal K}(P))$ is
\begin{equation*}
\{\Pi_{i=1}^t (n_i+1)!\}\Pi_{u=1}^{\infty} m_u!, 
\end{equation*}
where $n_i=|P_i|$ and $m_u$ is the number of parts of size $u$.

\end{proposition}
\proof We use again decomposition given in Theorem \ref{ProductDecomposition}. The symmetry group of the simplex of dimension $n$ has size $(n+1)!$, which yields the first term of the product. The second term comes from possible interchange of elements, if size of components are equal. \qed

Given a simplicial complex ${\cal K}$ of type $\{3,4\}$ and a $n$-face $\Delta=\{1,\dots, n+1\}$ of ${\cal K}$, 
define $g_i$ to be the reflection along the $(n-1)$-face $V_{n+1}-\{i\}$. The group, generated by $g_i$, is independent on the $n$-face $\Delta$ and it is denoted $Cox({\cal K})$. 

\begin{proposition}
Given a simplicial complex ${\cal K}={\cal K}(P)$ of type $\{3, 4\}$, the group $Cox({\cal K})$ has the following properties:

(i) $Cox({\cal K})$ is a Coxeter group with the following defining relations:
\begin{equation*}
\begin{array}{rcl}
g_i^2=1,\mbox{~~~~~~~}(g_ig_j)^2&=&1\mbox{~if~}i'\not= j',\\
                     (g_ig_j)^3&=&1\mbox{~if~}i'= j'.
\end{array}
\end{equation*}

(ii) $Cox({\cal K})$ is isomorphic to $\Pi_{i=1}^t Sym(1+|P_i|)$; its order is $\Pi_{i=1}^t (n_i+1)!$. 

(iii) $Cox({\cal K})$ is equal to $Aut({\cal K})$ if and only if all parts of the partition have different size (this case includes simplex and the bipyramid on a simplex).

(iv) $Cox({\cal K})$ is transitive on facets. 

(v) The fundamental domain of $Cox({\cal K})$ is a face if and only if ${\cal K}$ is the hyperoctahedron (i.e., the action of $Cox({\cal K})$ is regular on the $n$-faces).
If the complex ${\cal K}$ is different from hyperoctahedron, then $Cox({\cal K})$ is at least two times transitive on facets.
In general, the fundamental domain of $Cox({\cal K})$ is a simplex with angles $\frac{\pi}{q}$ for $q=2$ or $3$.
\end{proposition}

Clearly, if ${\cal K}$ is $(n+1)$-simplex, then $Cox({\cal K})$ is the irreducible group $A_n$. For all other simplicial complexes of type $\{3,4\}$, this group is an {\em reducible} Coxeter group.

The first case, when $Aut({\cal K})$ is not generated by reflections, appears for the complex ${\cal K}(\{1,2\}, \{3,4\})$. In general, if $Aut({\cal K})$ of a complex of type $\{3,4\}$ is generated by reflections, then it is a Coxeter group.
%(NEED TO CHECK ABOVE STATEMENT)

\begin{proposition}
Let ${\cal K}$ and ${\cal K'}$ be two simplicial complexes of type $\{3,4\}$, such that $Cox({\cal K})$ is isomorphic to $Cox({\cal K'})$. Then ${\cal K}$ and ${\cal K'}$ are isomorphic.
\end{proposition}
\proof We express ${\cal K}$ (respectively, ${\cal K'}$) as 
${\cal K}(P)$ (respectively, ${\cal K}(P')$) and denote by $m_u$ (respectively, by $m'_u$) the number of parts in $P$ (respectively, $P'$) of size $i$.

The group $Cox({\cal K})$ is isomorphic to $(Sym(2))^{m_1}\times (Sym(k))^{m_{k-1}}\times \dots$. The group $A_{k}$ is simple if $k\geq 5$; its multiplicity in the Jordan-H\"older decomposition of $Cox({\cal K})$ is $m'_{k-1}$. So, one has $m_k=m'_k$ if $k\geq 4$. The multiplicities of $C_2$ (respectively, $C_3$) in the decomposition of $Cox({\cal K})$ is $m_1+m_2+3m_3+\sum_{k\geq 4} m_k$ (respectively, $m_2+m_3$).
One has trivially $n+1=\sum_{k} km_k$. So, by solving the linear system, one obtains $m_k=m'_k$ if $k\geq 1$. \qed

A polytope is called {\em regular-faced} if all its facets are regular polytopes.
\begin{proposition}\label{TheoremBlBl}
Amongst simplicial complexes of type $\{3,4\}$, the only ones admitting regular-faced polytopal realization are two regular ones (simplex and hyperoctahedron) and bipyramid over simplex.
\end{proposition}
\proof First, we remind that any simplicial complex of type $\{3,4\}$ admits polytopal realization as a convex polytope.
All regular-faced polyhedra are known. All $92$ $3$-dimensional ones are classified in \cite{Jo1}.
All ones of higher dimension are classified in \cite{BlBl2}, \cite{BlBl} and references 1, 2 therein.

Besides two infinite families (pyramid over hyperoctahedron and bipyramid over simplex), the list of regular-faced, but not regular, polytopes given in \cite{BlBl}, contains only polytopes in dimension $4$. For dimension $3$ and $4$, it is easy to check the proposition. \qed

\begin{conjecture}
The number of simplicial complexes {\em of type $\{3,4,5\}$} (i.e., such that the link of any $(n-2)$-face is $C_3$, $C_4$ or $C_5$) is finite in any fixed dimension $n$.

\end{conjecture}
Remark that simplicial complexes of type $\{3,4,5\}$ can be non-spherical. For example, take the icosahedron and identify opposite vertices, edges and faces. The obtained complex has type $\{5\}$, skeleton $K_6$ and can be realized into the projective plane.

\begin{table}
\begin{center}
\begin{tabular}{||c||c|c||c|c||c||}
\hline
Partition $P$ with             &skeleton      &\# facets&$|Aut({\cal K})|$&\# orbits&$|Cox({\cal K})|$\\
${\cal K}={\cal K}(P)$         &$G({\cal K})$ &         &       &on vertices  &\\
\hline
$\{1,2,3\}$                    &$K_4$         &$4$      &$24$        &$1$       &$24$\\
$\{1\},\{2,3\}$                &$K_5-K_2$     &$6$      &$12$        &$2$       &$12$\\
$\{1\},\{2\},\{3\}$            &$K_6-3K_2$    &$8$      &$48$        &$1$       &$8$\\
\hline
$\{1,2,3,4\}$                  &$K_5$         &$5$      &$120$       &$1$       &$120$\\
$\{1\},\{2,3,4\}$              &$K_6-K_2$     &$8$      &$48$        &$2$       &$48$\\
$\{1,2\},\{3,4\}$              &$K_6$         &$9$      &$72^*$      &$1$       &$36$\\
$\{1\},\{2\},\{3,4\}$          &$K_7-2K_2$    &$12$     &$48$        &$2$       &$24$\\
$\{1\},\{2\},\{3\},\{4\}$      &$K_8-4K_2$    &$16$     &$384$       &$1$       &$16$\\
\hline
$\{1,2,3,4,5\}$                &$K_6$         &$6$      &$720$       &$1$       &$720$\\
$\{1\},\{2,3,4,5\}$            &$K_7-K_2$     &$10$     &$240$       &$2$       &$240$\\
$\{1,2\},\{3,4,5\}$            &$K_7$         &$12$     &$144$       &$2$       &$144$\\
$\{1\},\{2\},\{3,4,5\}$        &$K_8-2K_2$    &$16$     &$192$       &$2$       &$96$\\
$\{1\},\{2,3\},\{4,5\}$        &$K_8-K_2$     &$18$     &$144^*$     &$2$       &$72$\\
$\{1\},\{2\},\{3\},\{4,5\}$    &$K_9-3K_2$    &$24$     &$288$       &$2$       &$48$\\
$\{1\},\{2\},\{3\},\{4\},\{5\}$&$K_{10}-5K_2$ &$32$     &$3840$      &$1$       &$32$\\
\hline
\end{tabular}
\end{center}
\caption{All simplicial complexes ${\cal K}$ of type $\{3,4\}$ of dimension at most $4$}
\label{ListOfFirstCases}
\end{table}

In Table \ref{ListOfFirstCases}, we give details for simplicial complexes of type $\{3,4\}$ of small dimension. In this table we mark by $*$ the cases, where the group is not Coxeter. The orbits of vertices are computed with respect to the group $Aut$.

Two different simplicial complexes of type $\{3,4\}$ with the same skeleton appear, starting from dimension $5$: ${\cal K}(\{1,2\}, \{3,4,5,6\})$ and ${\cal K}(\{1,2,3\}, \{4,5,6\})$ both have skeleton $K_8$.

\begin{remark}
The simplicial complex ${\cal K}(\{1,2\},\{3,4\})$ has the following properties:

(i) It cannot be realized as a convex polytope in $\R^4$, in such a way that each its facet is {\em regular} tetrahedron.
But this complex admits such embedding in $\R^5$.
Moreover, it embeds into $5$-simplex: in fact, into the simplicial complex formed by all $3$-dimensional faces of the $5$-simplex (apropos, the above simplicial complex is not a pseudomanifold).

(ii) It provides an example, that the theorem of Alexandrov (\cite{Al}) does not admits an analog in dimension $3$:

An abstract $n$-dimensional Euclidean simplicial complex is formed of simplexes and distances between vertices. If an abstract simplicial complex is realized as the complex formed by a set of points on the boundary of a polytope (i.e., a boundary complex), then it is homeomorphic to a $n$-sphere and the sum of angles at every vertex is lower or equal to the total angle of a sphere of dimension $n-1$ (i.e., it has {\em non-negative curvature}).

Alexandrov's theorem (\cite{Al}) asserts that any abstract Euclidean simplicial complex of dimension $2$, which is homeomorphic to a $2$-sphere and has non-negative curvature, can be realized in $\R^3$ as a boundary complex.

The complex ${\cal K}(\{1,2\},\{3,4\})$ has $6$ vertices. Let us put equal distance to all edges, i.e., assume that all facets are regular simplexes. It is easy to see that the obtained complex has non-negative curvature.
It can be realized in $\R^5$ by the regular $5$-simplex.
All $4$-dimensional polytopes, whose facets are regular $3$-simplexes, have been classified in \cite{BlBl2} (see, more generally, Proposition \ref{TheoremBlBl}) and ${\cal K}(\{1,2\},\{3,4\})$ is not one of them.

%In fact, we say that a simplicial complex, formed by Euclidean
%simplexes is of {\em non-negative curvature} if the volume of 
%$\epsilon$-neighborhood of any point is less or equal than its
%Euclidean volume.
%is an example of a $3$-dimensional simplicial, which is homeomorphic
%to the $3$-sphere and has non-negative curvature but not admit such
%convex realization in $\R^4$.

\end{remark}

\section{Cubical complexes}
A {\em cubical complex} is a lattice, whose facets are combinatorial hypercubes. So, all its proper faces are combinatorial hypercubes too.

We are interested, especially, by cubical complexes of type $\{3,4\}$.
%Let us a call a cubical complex of type $\{3,4\}$ if the link of every $(n-2)$-face has length $3$ or $4$.

The hypercubes are only cubical complexes, such that any $(n-2)$-face belongs exactly to three $n$-faces.

The {\em star} of a vertex in a given complex is the subcomplex formed by all faces, which are incident to a given vertex.
The star of any vertex of a cubical complex of type $\{3,4\}$ is a {\em simplicial} complex of type $\{3,4\}$; so, the classification of such complexes in Corollary \ref{ConsequenceFundamentalLemma} characterize them also, but only locally.

\begin{proposition}\label{Case4CubicalComplex}
Let ${\cal K}$ be a cubical complex, such that the link of every $(n-2)$-face has size $4$; then ${\cal K}$ is non-spherical and, moreover:

(i) if ${\cal K}$ is simply-connected, then it is the cubical lattice $Z^n$;

(ii) otherwise, ${\cal K}$ can be obtained as a quotient of $Z^n$ by a torsion-free (i.e., without fixed points) subgroup of the symmetry group of $Z^n$ (i.e., the semidirect product of the Coxeter group $B_n$ and of the group of translations).
\end{proposition}
\proof If ${\cal K}$ is a cubical complex, whose $(n-2)$-faces are contained in exactly four $n$-faces, then the star of any vertex is $(n-1)$-hyperoctahedron; so, one has an unique way to extend it locally to a cubical complex. 
The simple-connectedness ensures that this construction will not repeat itself. Therefore, one gets the cubical lattice.

If ${\cal K}$ is a cubical complex, then its universal cover is the cubical lattice $Z^n$. So, ${\cal K}$ is obtained as the quotient of $Z^n$ by a torsion-free (i.e. having no fixed points) subgroup of $Aut(Z^n)$.
The sphere is simply-connected; so, it cannot be obtained as a proper quotient. \qed

Proposition \ref{Case4CubicalComplex}.(ii) gives, for example, cubical complex of type $\{4\}$ on torus and Klein bottle.

%\begin{problem}
%Characterize all polytopal cubical complexes of type $\{3,4\}$.
%\end{problem}

A $2$-dimensional cubical complex is called {\em quadrillage}. The polytopal quadrillages of type $\{3,4\}$ are exactly dual {\em octahedrites}, studied in \cite{DS} and \cite{DDS}, i.e., finite quadrangulations, such that each vertex has valency $3$ or $4$. So, the number of such complexes is not finite already in dimension two.
All dual octahedrites, which are isohedral, are: Cube, dual Cuboctahedron and dual Rhombicuboctahedron.

\begin{definition}
Let ${\cal K}$ be a cubical complex; define a {\em zone} as a circuit of $(n-1)$-faces of ${\cal K}$, where any two consecutive elements are opposite faces of a facet.

%(ii) the $p$-inflation of a zone is defined as the cubical complex obtained by splitting every $(n-1)$-face $F$ of the zone in $p^{n-1}$ hypercube in square fashion. NEED TO ELIMINATE
%(iii) a {\em railroad} is a circuit of $(n-2)$-faces contained in exactly $n$-faces and such that any $(n-2)$-face is adjacent to its neighbors on its opposite $(n-3)$-faces.

\end{definition}

The notion of zone corresponds, in the case of octahedrites, to the notion of central circuits (see \cite{DS} and \cite{DDS}).

Each zone of an $n$-dimensional hypercube corresponds to a pair of opposite facets (i.e., to a parallel class of edges). By cutting those edges by $m$ equispaced parallel planes and doing so for each zone, one obtains another example of a cubical polytopal complex.

%The $n$-hypercube has $n$ zones; we can split 
%in fact, the inflation is a polytopal cubical complex, since we can express this inflation by breaking each edge in $p$ elements and moving them by a small distance.
%Remark that, in some cases (for example, $n$-hypercube), we can cut the $n-1$ edges of the $(n-1)$-face independently, i.e., in $p_1$, \dots, $p_{n-1}$ edges. ?RECOVER CONSTRUCTION OF MISHA OF A POLYTOPAL COMPLEX
%Call a cubical complex {\em tight} if it is not inflation of any other cubical complex.

\section{Embeddability of skeletons of complexes in hypercubes}

\begin{theorem}\label{OldTheorem}
Let ${\cal K}$ be a closed simplicial complex of dimension $n\geq 3$.
Then one has:

(i) the skeleton of ${\cal K}$ is not embeddable, if ${\cal K}$
has an $(n-2)$-face belonging to at least five $n$-simplexes and such that
its link is an isometric cycle in the skeleton;

(ii) the skeleton of ${\cal K}$ is embeddable if ${\cal K}$ is of type $\{3,4\}$.
 \end{theorem}
\proof If there exists a $(n-2)$-face, such that its link has size at least
 six, then the skeleton ${\cal K}$ is not 5-gonal, since it contains
the isometric subgraph $K_5-K_3$.
If a $(n-2)$-face has a link of size five, then the skeleton of ${\cal K}$
contains the isometric subgraph $K_7-C_5$, which is not embeddable. 

All skeletons of simplicial complexes of type $\{3,4\}$ are of the form $K_m-hK_2$ and so, (\cite{DL}, Chapter 7.4) embeddable. \qed

For example, Theorem \ref{OldTheorem}.(i) implies non-embeddability
of following $3$-dimensional simplicial complexes:

(a) regular $4$-polytope 600-cell, since the link of each of its edge
is an isometric $C_5$;

(b) the skeleton of Delaunay partition of the {\em body centered cubic lattice}
(denoted also $A_3^*$), since the link of some edges is an isometric
$C_6$ (this skeleton is, moreover, not $5$-gonal).

\begin{figure}
\centering
\begin{minipage}{6cm}
\begin{equation*}
\begin{array}{ccc}
1267&  2346&  1246\\
2367&  1456&  1357\\
3467&  1237&  1234\\
4567&  3457&  1345\\
5167&      &
\end{array}
\end{equation*}\par
Facet-set of the simplicial complex ${\cal K}$
\end{minipage}
\begin{minipage}{6cm}
\epsfxsize=46mm
\epsffile{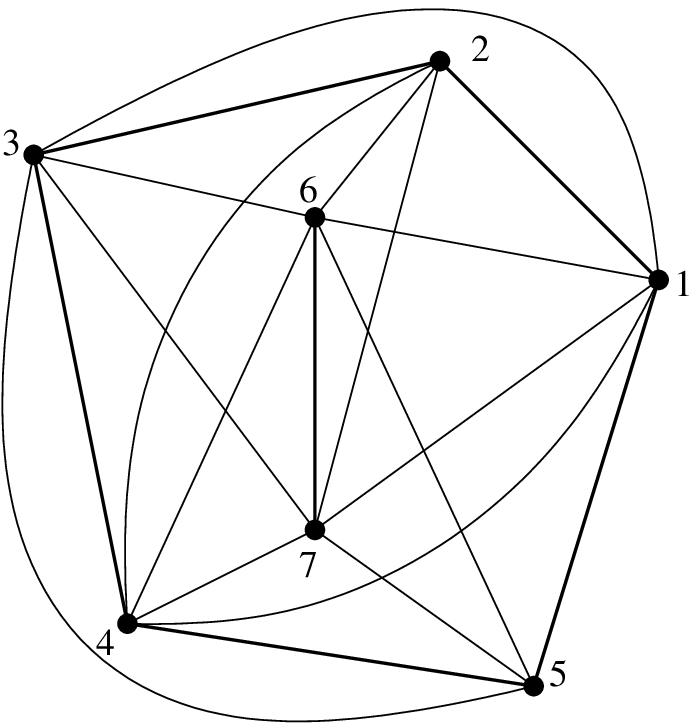}
\end{minipage}
\caption{Embeddable triangulation ${\cal K}$ with skeleton $K_7-K_2$ and non-isometric link $C_5=(1,2,3,4,5)$ of the edge $(6,7)$}
\label{CounterExample}
\end{figure}

The condition of isometricity of the link in Theorem \ref{OldTheorem}
is necessary (it was missed in \cite{DS3}). For example, there exists
an {\em embeddable} $3$-dimensional simplicial complex having an edge,
which belongs to five tetrahedra; its skeleton is $K_7-K_2$
(see Figure \ref{CounterExample}). The same graph $K_7-K_2$
appears also as the skeleton of a $n$-dimensional simplicial
complex of type $\{3,4\}$, but only for $n=4$.
 
\begin{corollary}
The skeleton of a finite cubical complex is embeddable and, moreover, with scale $1$, if and only if its path-metric satisfies the $5$-gonal inequalities.
\end{corollary}
\proof Clearly, such skeletons are bipartite graphs, and so, the embeddability of the skeleton implies that it is an isometric subgraph of some hypercube or, if infinite, of some cubic lattice $Z^m$. The result then follows from the characterization of isometric subgraphs of hypercubes, obtained in \cite{Djo} and reformulated in \cite{Avis}. \qed

On embeddings of {\em quadrillages} (i.e., cubical complexes of dimension two), we can say more.

One can check that any plane bipartite graph (i.e., all face-sizes are even) is an isometric subgraph of an hypercube if and only if all zones are {\em simple} (i.e., have no self-intersection) and each of them is convex (i.e., any two of its vertices are connected by a shortest path belonging to the zone).
So, the skeleton of a quadrillage is embeddable if and only if its zones are convex (and so, simple).

\begin{figure}
\centering
\begin{minipage}[b]{4.6cm}%
\centering
\epsfig{figure=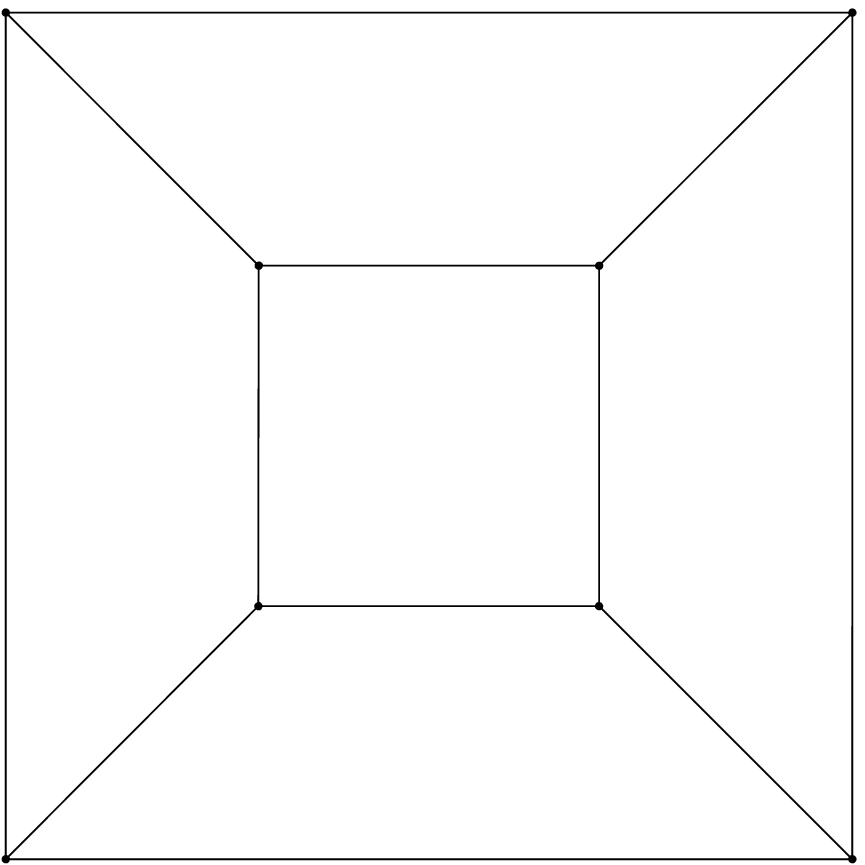,height=4cm}
$G_0$
\end{minipage}
\begin{minipage}[b]{4.6cm}%
\centering
\epsfig{figure=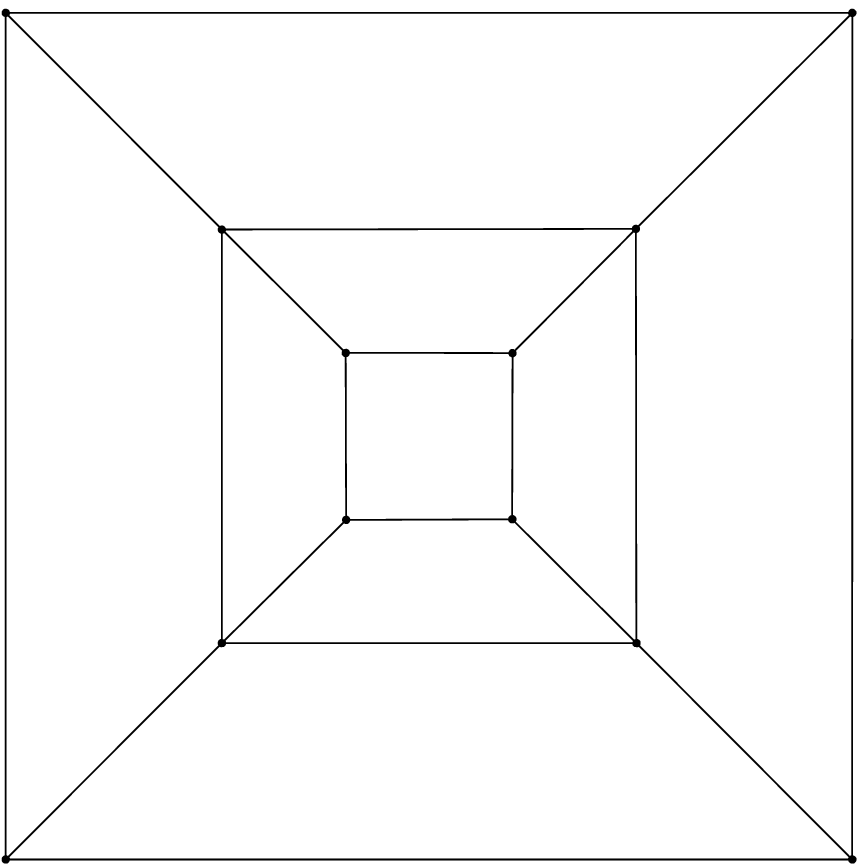,height=4cm}
$G_1$
\end{minipage}
\begin{minipage}[b]{4.6cm}%
\centering
\epsfig{figure=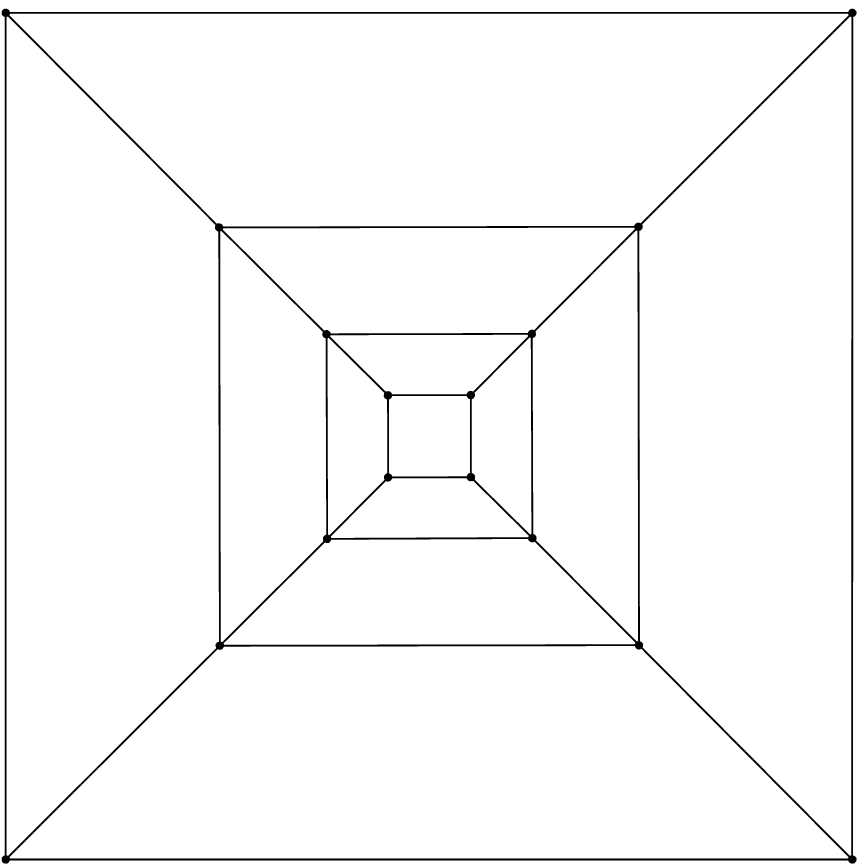,height=4cm}
$G_2$
\end{minipage}
\caption{The graphs $G_{t}$ (embeddable into $H_{t+3}$) for $t=0,1,2$}
\label{SeriesOfEmbeddables}
\end{figure}

All known embeddable polyhedral quadrillages (i.e., dual octahedrites, which are isometric subgraphs of hypercubes) are zonohedra:
dual Cuboctahedron and the family $G_t$ (for any integer $t\geq 0$), illustrated on the Figure \ref{SeriesOfEmbeddables} for the cases $t=0,1,2$. Graphs $G_t$ are embeddable in $H_{t+3}$.

% $(2-Prism^t_4)^*$ (i.e., dual $t$-inflation of the octahedron on one of its central circuits or cutting of railroad of the Cube, see drawing below NEED TO BE DONE); 

\begin{figure}
\centering
\begin{minipage}[b]{8cm}%
\centering
\epsfig{figure=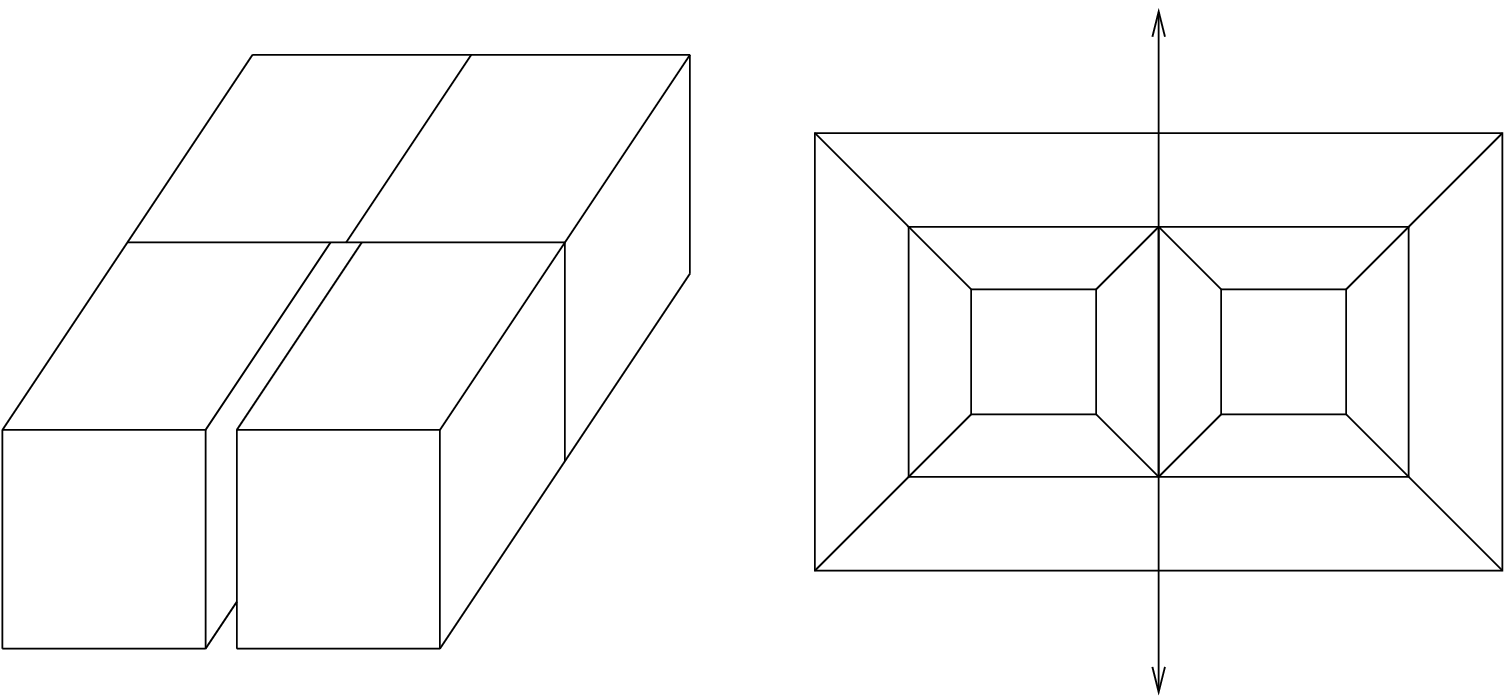,width=7cm}
\end{minipage}
\begin{minipage}[b]{8cm}%
\centering
\epsfig{figure=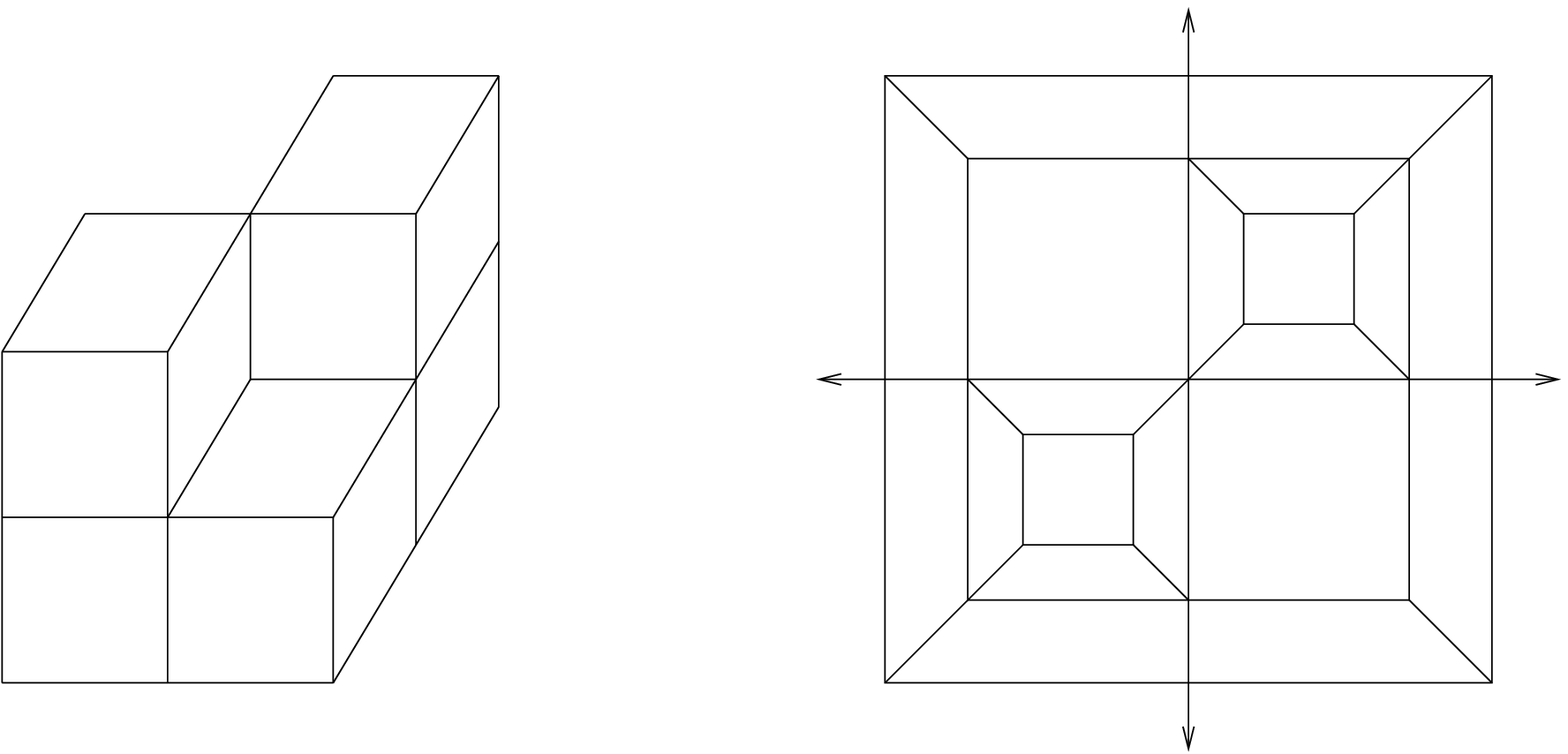,width=7cm}
\end{minipage}
\caption{Two non-embeddable quadrillages with simple zones}
\label{CntExamples}
\end{figure}

%A quadrillage is {\em topologically embeddable} if it is topologically embeddable in a cubic lattice $Z^m$. 
The simplicity of all zones is a necessary condition for topological embedding of quadrillage in a cubical lattice $Z^m$; it is not sufficient even in the spherical case (see Figure \ref{CntExamples}). But a quadrillage, such that its skeleton is {\em isometric} subgraph of hypercube, is topologically embeddable in a cubical lattice.

\end{document}